\documentclass[12pt,english]{article}
\usepackage{times}
\usepackage[T1]{fontenc}
\usepackage[latin1]{inputenc}
\usepackage{graphicx}
\usepackage{amssymb}

\makeatletter

\newcommand{\noun}[1]{\textsc{#1}}
\providecommand{\tabularnewline}{\\}

\usepackage{graphicx}
\usepackage{amssymb}

\usepackage{graphicx}
\usepackage{amssymb}

\usepackage{graphicx}
\usepackage{amssymb}

\usepackage{graphicx}
\usepackage{amssymb}

\usepackage{graphicx}
\usepackage{amssymb}

\usepackage{graphicx}
\usepackage{amssymb}

\usepackage{array}
\usepackage{amssymb}

\usepackage{amsmath,amssymb,amscd,enumerate,amsthm}
\usepackage{graphics}
\usepackage{epsfig}

\newtheorem{theorem}{Theorem}[section]
\newtheorem{lemma}[theorem]{Lemma}

\newtheorem{algorithm}{Algorithm}[section]

\usepackage{babel}

\usepackage{babel}

\usepackage{babel}

\usepackage{babel}

\usepackage{babel}

\usepackage{babel}

\usepackage{babel}

\usepackage{babel}
\makeatother
\begin{document}

\title{A Sublinear Algorithm of Sparse Fourier Transform for Nonequispaced
Data\thanks{This
        work was partially supported by NSF grant DMS-03168875 and AFOSR grant 109-6047.}}

\author{Jing Zou %
\footnote{Program of Applied and Computational Mathematics, Princeton University,
Fine Hall, Washington Road, Princeton, NJ 08544, (\texttt{jzou@math.princeton.edu})%
}}
\date{}
\maketitle
\begin{abstract}
We present a sublinear randomized algorithm to compute a sparse Fourier
transform for nonequispaced data. Suppose a signal $S$ is known to
consist of $N$ equispaced samples, of which only $L<N$ are available.
If the ratio $p=L/N$ is not close to 1, the available data are typically
non-equispaced samples. Then our algorithm reconstructs a near-optimal
$B$-term representation $R$ with high probability $1-\delta$, in
time and space $poly(B,\log(L),\log p,\log(1/\delta),$ $\epsilon^{-1})$,
such that $\Vert S-R\Vert^{2}\leq(1+\epsilon)\Vert S-R_{opt}^{B}\Vert^{2}$,
where $R_{opt}^{B}$ is the optimal $B$-term Fourier representation
of signal $S$. The sublinear $poly(\log L)$ time is compared to
the superlinear $O(N\log N+L)$ time requirement of the present best
known Inverse Nonequispaced Fast Fourier Transform (INFFT) algorithms.
Numerical experiments support the advantage in speed of our algorithm
over other methods for sparse signals: it already outperforms INFFT
for large but realistic size $N$ and works well even in the situation
of a large percentage of missing data and in the presence of noise. 
\end{abstract}

\section{Introduction}

We consider the problem in which the recovery of a discrete time signal
$S$ of length $N$ is sought when only $L$ signal values are known.
In general, this is of course an insoluble problem; we consider it
here under the additional assumption that the signal has a sparse
Fourier transform. Let us fix the notations: the signal is denoted
by $S=(S(t))_{t=0,\ldots,N-1}$, but we have at our disposal only
the $(S(i))_{i\in T}$, where the set $T$ is a subset of $\{0,\ldots,N-1\}$
and $|T|=L$. The Fourier transform of signal $S$ is $\hat{S}=(\hat{S}(0),\ldots,\hat{S}(N-1))$,
defined by $\hat{S}(\omega)=\frac{1}{\sqrt{N}}\sum_{t=0}^{N-1}S(t)e^{-2\pi i\omega t/N}$.
In terms of the Fourier basis functions $\phi_{\omega}(t)=\frac{1}{\sqrt{N}}e^{2\pi i\omega t/N}$,
$S$ can be written as $S=\sum_{\omega=0}^{N-1}\hat{S}(\omega)\phi_{\omega}(t)$;
this is the (discrete) Fourier representation of $S$. A signal $S$
is said to have a $B$-sparse Fourier representation, if there exists
a subset $\Omega\subset\{0,\ldots,N-1\}$ with $|\Omega|=B$, and
values $c(\omega)\neq0$ for $\omega\in\Gamma$, such that $S(t)=\sum_{\omega\in\Omega}c(\omega)\phi_{\omega}$.
For a signal that does not have a $B$-sparse Fourier representation,
we denote by $R_{opt}^{B}(S)$ the optimal $B$-term Sparse Fourier
representation of $S$.

This paper presents a sublinear algorithm to recover a $B$-sparse
Fourier representation of a signal $S$ from incomplete data. Our
algorithm also extends to the case where the Fourier transform $\hat{S}$
is not $B$-sparse, where we aim to find a near-optimal $B$-term
Fourier representation, i.e. $R=\sum_{\omega\in\Gamma}c(\omega)\phi_{\omega}$,
such that \begin{equation}
\| S-R\|=\Vert S-\sum_{\omega\in\Gamma}c(\omega)\phi_{\omega}\Vert_{2}^{2}\leq(1+\epsilon)\Vert S-R_{opt}^{B}(S)\Vert_{2}^{2}.\end{equation}

A typical situation where our study applies is the observation of
non-equispaced data, where the samples are nevertheless all elements
of $\tau\mathbb{Z}$ for some $\tau>0$. For a signal with evenly
spaced data, the famous Fast Fourier Transform (FFT) computes all
the Fourier coefficients in time $O(N\log N)$. However, the requirement
of equally distributed data by FFT raises challenges for many important
applications. For instance, because of the occurrence of instrumental
drop-outs, the data may be available only on a set of non-consecutive
integers. Another example occurs in astronomy, where the observers
cannot completely control the availability of observational data:
a telescope can only see the universe on nights when skies are not
cloudy. In fact, computing the Fourier representation from irregularly
spaced data has wide applications \cite{Ware} in processing astrophysical
and seismic data, the spectral method on adaptive grids, the tracking
of Lagrangian particles, and the implementation of semi-Lagrangian
methods.

In many of these applications, a few large Fourier coefficients already
capture the major time-invariant wave-like information of the signal,
and we can thus ignore very small Fourier coefficients. To find a
small set of the largest Fourier coefficients and hence a (near) optimal
$B$-sparse Fourier representation of a signal that describes most
of the signal characteristics is a fundamental task in applied Fourier
Analysis.

An equivalent version of this problem is as follows: define the matrix
$A:=(e^{2\pi ikt_{j}})_{k=0,\ldots,N;}$ $_{j=0\ldots,L-1}$, where
the $t_{j}$ are the locations of the available samples. Given $S(t_{j})$,
we want to reconstruct the signal $S$, or equivalently, its Fourier
coefficients $\hat{S}_{k}$, so that $A\hat{S}=S$. This linear system
is over-determined. Several algorithms \cite{Bjork}\cite{Hanke}
\cite{KP} have provided efficient approaches to solve this problem.
Among all INFFT algorithms, the iterative CGNE approach of \cite{FGS}
in the benchmark software NFFT 2.0 is one of the fastest methods;
it takes time $O(L^{1+(d-1)/\beta}\log L)$, where $L$ is the number
of available points, $d$ is the number of dimensions, and $\beta>1$
is the smoothness for the original signal. The super-linearity relationship
between the running time and $N$ (recall $L=pN$, where $p$ is the
percentage of available data) poses difficulties in processing large
dimensional signals, which have nothing to do with the unequal spacing.
It follows that identifying a sparse number of significant modes and
amplitudes is expensive for even fairly modest $N$. Our goal in this
paper is to discuss much faster (sublinear) algorithms that can identify
the sparse representation or approximation with coefficients $a_{1},\ldots,a_{B}$
and modes $\omega_{1},\ldots,\omega_{B}$ for unevenly spaced data.
These algorithms will not use all the samples $S(0),\ldots,S(N-1)$,
but only a very sparse subset of them.

Our approach is based on the paper \cite{GGIMS} that shows how to
construct the Fourier representation for a signal $S$ with $B$-sparse
Fourier representation in time and space $poly(B,\log N,$ $1/\epsilon,\log(1/\delta))$
on equal spacing data. The algorithm contains some random elements
(which do not depend on the signal); their approach guarantees that
the error of estimation is of order $\epsilon\Vert S\Vert^{2}$ with
probability exceeding $1-\delta$. The ideas in \cite{GGIMS} have
also been applied by its authors to sparse wavelet, wavelet packet
representation, and histograms \cite{GGIKMS}. We have dubbed the
whole family of algorithms RA$\ell$STA (for Randomized Algorithm
for Sparse Transform Approximation); when dealing only with Fourier
Transforms, as is the case here, we specialize it to RA$\ell$SFA
(F for Fourier). Zou, Gilbert, Strauss and Daubechies \cite{Zou}
improved and implemented the algorithm greatly. It convincingly beats
FFT when the number of grid points $N$ is reasonably large. The crossover
point lies at $N\simeq25,000$ in one dimension, and at $N\simeq460$
for data on a $N\times N$ grid in two dimensions for a two-mode signal.
When $B=13$, RA$\ell$SFA surpasses $FFT$ at $N\geq300,000$ for
one dimensional signals and $1100$ for two dimensional signals.

In this paper, we modify RA$\ell$SFA to solve the irregularly spaced
data problem. The new NERA$\ell$SFA (Nonequispaced RA$\ell$SFA)
uses sublinear time and space $poly(B,\log L,\epsilon,\log(1/\delta),$
$\log p)$ to find a near-optimal $B$-term Fourier representation,
such that $\Vert S-R\Vert^{2}\leq(1+\epsilon)\Vert S-R_{opt}\Vert^{2}$
with high probability $1-\delta$. Similar to the RA$\ell$SFA algorithm,
it outperforms existing INFFT algorithms in processing sparse signals
of large size.

\textbf{Notation and Terminology} Denote by $\chi_{T}$ a signal that
equals 1 on a set $T$ and zero elsewhere in the time domain. We say
a signal $H$ is $q$ percent pure, if there exists a frequency $\omega$
and a signal $\rho$, such that $H=ae^{2\pi i\omega t/N}+\rho$, with
$|a|^{2}\geq(q\%)\| H\|^{2}$. To quantify the unevenness of the data,
introduce a parameter $p=L/N$ to be the percentage of the available
data over all the data, where $L$ is the number of available data.
Obviously a larger $p$ corresponds to more information about the
signal. We use $L^{2}$-norm throughout the paper, which is denoted
by $\|.\|$. The convolution $F*G$ is defined as $F*G(t)=\sum_{s}F(s)G(t-s)$.
It follows that $\widehat{F*G}(\omega)=\sqrt{N}\hat{F}(\omega)\hat{G}(\omega)$. 

A Box-car filter with width $2k+1$ is defined as follows: \begin{eqnarray*}
\chi_{k}(t) & = & \left\{ \begin{array}{cc}
\frac{\sqrt{N}}{2k+1} & \,\,\,\,\,\,\,\textrm{$if\,\,\,\,\,-k\leq t\leq k$ },\\
0 & \,\,\,\,\,\,\,\,\,\, if\,\,\, t>k\,\, or\,\, t<-k\end{array}\right.\end{eqnarray*}
 In the frequency domain, this filter is in the form of \begin{equation}
\hat{\chi}_{k}(\omega)=\left\{ \begin{array}{c}
\frac{sin((2k+1)\pi\omega/N)}{(2k+1)sin(\pi\omega/N)}\,\,\,\,\,\,\,\,\, if\,\,\omega\neq0\\
\,\,\,\,\,\,\,\,\,\,1\,\,\,\,\,\,\,\,\,\,\,\,\,\,\,\,\,\,\,\,\,\,\,\,\,\,\, if\,\,\,\omega=0\end{array}\right.\end{equation}

A dilation operation on signal $H$ with a dilation factor $\sigma$
is defined as $H^{(\sigma)}(t)=H(\sigma t)$ for every points $t$.

\textbf{Organization} The paper is organized as follows. In Section
2, we give the outline of the RA$\ell$SFA algorithm. Section 3 presents
the modification of RA$\ell$SFA that deals with the unavailability
of some samples by a greedy method. In Section 4, an interpolation
technique is introduced for better performance. Finally, we compare
numerical results with existing algorithms in Section 5.

\section{Set-up of RA$\ell$SFA}

Given a signal $S$ of length $N$, the optimal $B$-term Fourier
representation $R_{opt}^{B}(S)$ uses only $B$ frequencies; it is
simply a truncated version of the Fourier representation of $S$,
retaining only the $B$ largest coefficients. The following theorem
is the main result of \cite{GGIMS}.

\begin{theorem}
Let an accuracy factor \(\epsilon\), a failure probability \(\delta\), and a sparsity target \(B \in \mathbb{N}, B \ll N\)
be given. Then for an arbitrary signal \(S\) of length \(N\), RA\(\ell\)SFA will find a
\(B\)-term approximation \(R\) to \(S\), at a cost in time and space of
order \(poly(B,\log(N), 1/\epsilon, \log(1/\delta) )\) and with probability exceeding \(1-\delta\), so that \(\|S-R\|^2 \leq (1+\epsilon)\|S-R_{opt}^B(S)\|^2_2\).
\end{theorem}\addtocounter{algorithm}{+1}

The striking fact is that RA$\ell$SFA can build a near-optimal representation
$R$ in sublinear time $poly(\log N)$ instead of the $O(N\log N)$
time requirement of other algorithms. Its speed surpasses FFT as long
as the length of a signal is sufficiently large. If a signal is composed
of only $B$ modes, RA$\ell$SFA constructs $S$ without any error.

The main procedure is a Greedy Pursuit with the following steps:

\begin{algorithm} \label{alg:total1}{\textsc{Total Scheme \cite{Zou}}}
\begin{enumerate}
    \item Initialize the representation signal \(R\) to 0. Set the maximum number of iterations \(ITER=B\log(N)\log(1/\delta)/\epsilon^{2}\).
    \item Test whether \(\Vert S-R \Vert\) appears to be less than some user threshold, $\iota$. 
      If yes, return the representation signal \(R\) and the whole algorithm ends; else go to step 3..
    \item Locate Fourier Modes \(\omega\) for the signal \(S-R\) by isolation and group test procedures.
    \item Estimate Fourier Coefficients at \(\omega\): \(\widehat{(S-R)}(\omega)\).
    \item Update the representation signal \(R\leftarrow R+\widehat{(S-R)}(\omega) \phi_{\omega}(t)\).
    \item If the total number of iterations is less than \(ITER\), go to 2; else return the representation \(R\).
  \end{enumerate}\label{alg:total1}

\end{algorithm}\addtocounter{theorem}{+1}

The basic idea of Algorithm \ref{alg:total1} is to identify significant
frequencies and then estimate their corresponding coefficients. In
order to locate those nonzero frequencies, we first construct a new
signal where a previous significant frequency becomes predominant.
Then a recursive approach called group test finds the exact label
of this predominant mode, by splitting intervals, comparing energies,
and keeping only intervals with large energies. After the frequency
is located, coefficient estimation procedures give a good estimation
by taking means and medians of random samples.

\section{NERA$\ell$SFA with Greedy Technique}

RA$\ell$SFA samples from a signal, implicitly assuming that uniform
and random sampling is possible, with a fixed cost per sample. This
raises challenges for processing unevenly spaced data. Specifically
speaking, Fourier coefficients and norms can not be estimated properly.
Thus one has to modify steps 3 and 4 accordingly. In this section,
NERA$\ell$SFA, a modified version of RA$\ell$SFA with greedy technique,
is introduced to overcome these problems.

The basic idea is a greedy pursuit for an available data point. Whenever
the algorithm samples at a missing data point, it searches some other
random indices $t$ until it finds one available data point $S(t)$
as the substitute. This technique is used in estimating both Fourier
coefficients and norms.

A good data structure is important to save running time cost. We denote
the availability of a data point by a label, say +1 for available
and 0 for unavailable. Hence, the label is tested to see if its corresponding
sample is valid. An alternative solution is to store all the sorted
labels of available data in a long list. However, each search takes
time $O(\log(N))$, which introduces a $O(\log N)^{2}$ factor into
the whole computation. As the empirical results show, the running
time of NERA$\ell$SFA algorithm is linear to $\log N$. For this
reason, we selected the first method.

We now give a more detailed discussion of the different procedures
used in steps 3 and 4 of Algorithm \ref{alg:total1}.

\subsection{Estimating Fourier Coefficients}

First, we give the procedure for estimating Fourier coefficients for
unevenly spaced data as follows.

\begin{algorithm} {\textsc{Estimating Individual Fourier Coefficients}} \label{estcoef} \\
Input a signal \(S\), a frequency \(\omega\), \(n=2\log(1/\delta)\), \(m=8/\epsilon^2\).
\begin{enumerate}
\item For \(i=1, \ldots, n\)
\item For \(j=1,\ldots, m\) \\
\text{}\hspace{10mm} Randomly generate the index \(t\) until \(S(t)\) is available. \\
\text{}\hspace{10mm} Then let \(t_{ij}=t\). Evaluate \(k(t_{ij})=<S(t_{ij}), \phi_{\omega}(t_{ij})>\).
\item Take the means of \(m\) samples \(k(t_{ij})\), i.e. \(p(i)=\sum_{j=1}^{m} {k(t_{ij})}\), where \(i=1,\ldots, n\).
\item Take the median of \(n\) samples \(c=median_i (p(i))\), where \(i=1,\ldots,n\).
\item Return \(c\) as the estimation of the Fourier coefficient \(\hat{S}(\omega)\).
\end{enumerate} \label{alg:coeff}
\end{algorithm}\addtocounter{theorem}{+1}

Next, we show that using unevenly spaced data leads to a very good
approximation to the true coefficient. The first lemma is one of most
fundamental theorems in randomized algorithms. It essentially states
that by repeating an experiment enough times, a small probability
event will happen eventually.

\begin{lemma} \label{lm:rept}
If an event happens with probability \(p\), then in the first \(k> \log \delta / \log(1-p)\) iterations, it happens at least once with success
probability \(1-\delta\).
\label{lm:rept}\end{lemma}\addtocounter{algorithm}{+1}

In our case, only $p=L/N$ percentage of the data is available, so
that $k>\log\delta/\log(1-L/N)$ trials are needed to generate one
available data point with success probability at least $1-\delta$.

In fact, most of the Fourier coefficients of a characteristic function
on a typical set $T$ are small, under some conditions. The following
lemma makes this more explicit.

\begin{lemma} \label{lm:smallfilter}
Suppose the components \(X_j\) of a discrete random variable \(X=(X_j)_{j=0} ^{N-1}\) are identically and independently distributed in \( \{0,1\} \), with \(p = Prob(X_j=1)\).
Define the random set \(T=\{j \in \{0, \ldots, N-1 \}  |X_j=1 \}\) to be the set of all available data; \(\hat{\chi}_T(\omega)\) is the Fourier transform of \(\chi_T(t)= \sum_{j=0}^{N-1}X_j\). If \(p \geq \frac{1}{1+(N-1) \lambda \tau^2}\),
then \\
    \begin{equation}
      Prob(|\hat{\chi}_T(\omega)|^2 \geq \lambda) \leq \tau^2.
      \end{equation}
\label{lm:smallfilter}\end{lemma}\addtocounter{algorithm}{+1} \begin{proof}
First, we claim that \(E(|\hat{\chi}_T(\omega)|^2) \leq \frac{(1-p)}{p(N-1)}\). \\
Since \(\hat{\chi}_T(\omega) = \frac{1}{pN}\sum_{j \in T}(e^{2 \pi i \omega j /N})\),
we have
\begin{eqnarray}
|\hat{\chi}_{T}(\omega)|^{2}=\frac{1}{p^{2}N^{2}}\sum_{j,k\in T}e^{2\pi i\omega(j-k)/N} \\
=\frac{1}{p^{2}N^{2}}\sum_{j\in T}1+\frac{1}{p^{2}N^{2}}\sum_{j,k\in T,j\neq k}e^{2\pi i\omega(j-k)/N}.
\nonumber \end{eqnarray}
It follows that
\begin{equation}
E(|\hat{\chi}_{T}(\omega)|^{2})=\frac{1}{pN}+\frac{1}{p^{2}N^{2}}p\frac{pN-1}{N-1}\sum_{j,k=0,j\neq k}^{N-1}e^{2\pi i\omega(j-k)/N}.
\nonumber \end{equation}
Observe that \(\sum_{j,k=0,j\neq k}^{N-1}e^{2\pi i\omega(j-k)/N}=|\sum_{j=0}^{N-1}e^{2\pi i\omega j/N}|^{2}-\sum_{j=0}^{N-1}1=(N\delta_{\omega,0})^{2}-N\), hence
\begin{eqnarray}
E(|\hat{\chi}_{T}(\omega)|^{2}) = \frac{1}{pN}+\frac{1}{pN^{2}}\frac{pN-1}{N-1}(N^{2}\delta_{\omega,0}-N)
=\frac{1}{pN}\left\{ 1+\frac{pN-1}{N-1}(N \delta_{\omega,0}-1)\right\}  \nonumber \\
=\frac{1}{pN(N-1)}\left\{ N-1+(pN-1)(N\delta_{\omega,0}-1)\right\}.  \nonumber
\end{eqnarray}
 By Markov's Inequality, when \(\omega \neq 0\), we have
 \begin{equation}
Prob(|\hat{\chi}_{T}(\omega)|^{2}\geq\lambda)\leq\frac{E(|\hat{\chi}_{T}(\omega)|^{2})}{\lambda}=\frac{1-p}{p(N-1)\lambda}.\nonumber \end{equation}
 Since \(p\geq\frac{1}{1+(N-1)\lambda \tau^{2}}\), it follows that
 \begin{equation}
Prob(|\hat{\chi}_{T}(\omega)|^{2}\geq\lambda)\leq\tau^{2}. \nonumber \end{equation}
 That is , for any \(\omega\neq 0\), with probability at least \(1-\tau^{2}\)
 \begin{equation}
|\hat{\chi}_{T}(\omega)|\leq \sqrt{\lambda}.
\end{equation}
\qquad \end{proof}

In particular, we want both $\lambda$ and $\tau$ to be small, meaning
that $p$ cannot be too small itself. 

Next, we consider the conditions for the two coefficients $\hat{S}(\omega)$
and $\hat{S}_{1}(\omega)=\widehat{S\cdot\chi_{T}}(\omega)$ to be
close.

\begin{lemma} \label{lm:diffest}
Suppose the parameters \(T\), \(S\), \(\chi_T(t)\), \(\lambda\), \(\tau\), \(p\) are as stated in Lemma \ref{lm:smallfilter}, and define \(S_1(t) = S(t)\chi_T(t)\). If \( p\geq \frac{1}{1+(N-1)\lambda \tau^2}\), and \(\tau \leq \sqrt{1-(1-\delta)^\frac{1}{B}}\),
then, for any \(\omega\),
\begin{equation}
|\hat{S}(\omega) - \hat{S}_1(\omega)| \leq \sqrt{B \lambda} \|S\|_2.
\end{equation}
with probability exceeding \(1-\delta\).
\label{lm:diffest}\end{lemma}\addtocounter{algorithm}{+1} \begin{proof}
Suppose the significant terms of signal \(S\) are \(\omega_i\), where \(i=1, \ldots, B\).\\
Since \(S_1(t) = S(t) \chi_T(t)\) and thus \(\hat{S}_1(\omega) =\hat{S}(\omega)* \hat{\chi}_T(\omega)\),
then
\begin{equation}
\begin{split}
\hat{S}_{1}(\omega_{j}) = \sum_{i=1}^{B}\hat{S}(\omega_{i})\hat{\chi}_T(\omega_{j}-\omega_{i})
=\hat{S}(\omega_{j})\hat{\chi}_{T}(0)+\sum_{i=1,\omega_{j}\neq\omega_{i}}^{B}\hat{S}(\omega_{i})\hat{\chi}_T(\omega_{j}-\omega_{i})  \\
=\hat{S}(\omega_{j})+\sum_{i=1,\omega_{j}\neq\omega_{i}}^{B}\hat{S}(\omega_{i})\hat{\chi}_T(\omega_{j}-\omega_{i}). \nonumber
\end{split}
\end{equation}
Therefore
\begin{equation}
|\hat{S}_{1}(\omega_{j})-\hat{S}(\omega_{j})| =  |\sum_{i=1,\omega_{j}\neq\omega_{i}}^{B}\hat{S}(\omega_{i})\hat{\chi}_T(\omega_{j}-\omega_{i})|
\end{equation}
\begin{equation}
\leq\sqrt{\sum_{i=1,\omega_{j}\neq\omega_{i}}^{B}|\hat{S}(\omega_{i})|^{2}}\sqrt{\sum_{i=1,\omega_{j}\neq\omega_{i}}^{B}|\hat{\chi}_{T}(\omega_{j}-\omega_{i})|^{2}} \nonumber \\
\leq\| S\|_2 \sqrt{\sum_{i=1,\omega_{j}\neq\omega_{i}}^{B}|\hat{\chi}_{T}(\omega_{j}-\omega_{i})|^{2}}.
\nonumber \end{equation}
 Because \(p\geq\frac{1}{1+(N-1)\lambda \tau^{2}}\), we have \(|\hat{\chi}_{T}(\omega)|^{2}\leq\lambda\)
with probability at least \(1-\tau^2\) for any \(\omega\neq0\). This implies that
 \(|\hat{S}_{1}(\omega_{j})-\hat{S}(\omega_{j})|\leq\| S\|_2 \sqrt{B \lambda}\) with probability at least  \((1-\tau^2)^B \geq (1-\delta)\) \\
Then
\begin{equation}
|\hat{S}_{1}(\omega_{j})-\hat{S}(\omega_{j})|\leq \sqrt{B \lambda} \| S\|_2.
\end{equation}
For those \(\omega\notin\{\omega_{i},i=1,\ldots,B\}\),
\begin{eqnarray}
\hat{S}_{1}(\omega)=\sum_{i=1}^{B}\hat{S}(\omega)\hat{\chi}_T(\omega-\omega_{i}), \nonumber \\
 \end{eqnarray}
and we conclude similarly that \(|\hat{S}_{1}(\omega)-\hat{S}(\omega)|\leq \sqrt{B \lambda} \| S\|_2.\), with probability at least \(1-\delta\).
\end{proof}

We shall use Algorithm \ref{alg:coeff} to estimate $\hat{S}_{1}(\omega)$;
we now look at how close the approximation $A$ (i.e. the output of
Algorithm \ref{alg:coeff}) of $\hat{S}_{1}(\omega)$ is to the true
coefficient $\hat{S}(\omega)$.

\begin{lemma}
\label{lm:coefftot} For a set of parameters \(T\), \(S\), \(\chi_T(t)\), \(\lambda\), \(\tau\), \(p\) as stated in Lemma \ref{lm:smallfilter}, if \( p\geq \frac{1}{1+(N-1)\lambda \tau^2}\), and \(\tau \leq \sqrt{1-(1-\delta)^{1/B} }\), then Algorithm \ref{alg:coeff} for signal \(S_1(t)=S(t) \* \chi_T(t)\) gives a good estimation \(A\) of \(\hat{S}(\omega)\), such that
 \begin{equation}
|A-\hat{S}(\omega)|\leq (\sqrt{\lambda} + \sqrt{B \lambda})\|S\|_2.
\end{equation}
 with high probability.\label{lm:coefftot}\end{lemma}\addtocounter{algorithm}{+1} \begin{proof}
Lemma 4.2 in \cite{Zou} says that the coefficient estimation algorithm returns \(A\), such that
\begin{equation}
|A-\hat{S}_1(\omega)|\leq \sqrt{\lambda} \|S\|_2.
\end{equation}
By Lemma \ref{lm:diffest}
\begin{equation}
|\hat{S}_{1}(\omega)-\hat{S}(\omega)|\leq \sqrt{B \lambda} \| S\|_2.
\end{equation}
Thus
\begin{equation}
|A-\hat{S}(\omega)| \leq |A-\hat{S}_1(\omega)|+|\hat{S}_{1}(\omega)-\hat{S}(\omega)|\leq(\sqrt{\lambda}+\sqrt{B \lambda})\| S\|_2.
\end{equation}
\end{proof}

Finally, we derive the conclusion about estimating coefficients.

\begin{theorem}
\label{lm:mycoeff2} For a set of parameters \(T\), \(S\), \(\chi_T(t)\), \(\lambda\), \(\tau\), \(p\) as stated in Lemma
\ref{lm:smallfilter}, if \(\lambda \leq \frac{\epsilon}{2(B+1)}\) and
\( p\geq \frac{1}{1+(N-1)\lambda \tau^2}\), then every application of Algorithm \ref{estcoef} produces, for each frequency \(\omega\) and each signal \(S\), and each \(\lambda>0\), with high
probability, an output \(A\) (after inputting \((S, \omega, \epsilon)\) ), such that \(|A-\hat{S}(\omega)|^2 \leq \epsilon \|S\|_2^2\).
\end{theorem}\addtocounter{algorithm}{+1} \begin{proof}
By Lemma \ref{lm:coefftot},
\begin{equation}
|A-\hat{S}(\omega)| \leq  (\sqrt{\lambda} + \sqrt{B \lambda}) \|S\|_2.
\end{equation}
Thus we have
\begin{equation}
|A-\hat{S}(\omega)|^2 \leq  2(\lambda + B \lambda) \|S\|_2^2.
\end{equation}
From the conditions \(2(\lambda + B \lambda) \leq  \epsilon \), it follows that
\begin{equation}
|A-\hat{S}(\omega)|^2 \leq \epsilon \|S\|_2^2.
\end{equation}
 \qquad\end{proof}

When we are able to get most of the data, the computational cost for
estimating Fourier coefficients on unevenly spaced data is only slightly
more than for the evenly spaced data case. The time to compute the
signal value remains almost the same as for the evenly spaced data
case. The extra time, in the worst case $O(\frac{\log\delta}{\epsilon_{1}^{2}p\log(1-p)})$,
comes from visiting unavailable data. Fortunately, the visit operation
is very fast and therefore contributes little to the total time, especially
when most of the data are available.

Moreover, as in \cite{Zou}, one can speed up the algorithm by using
multi-step coarse-to-fine coefficient estimation procedures, which
turns out to be more efficient than single-step accurate estimation;
the proof is entirely analogous to Lemma 4.3 in \cite{Zou}.

\subsection{Estimating Norms}

The basic idea for locating the label of a significant frequency is
to compare the energies (i.e. the $L^{2}$ norm) of signals restricted
in different frequency intervals. If the energy of some interval is
relatively large, the significant mode is in that region with higher
probability. We construct the following new signals to focus on certain
intervals \begin{equation}
H_{j}(t)=\chi_{1}(t)e^{\frac{2\pi ijt}{16}}\ast\chi_{[-q_{1},q_{1}]}(\sigma t)e^{\frac{2\pi it\theta}{N}}\ast S\end{equation}
 where 2$q_{1}+1$ is the filter width, $j=0,\ldots,15$, $\sigma$
and $\theta$ are random dilation and modulation factors. (Please
see \cite{Zou} for an explanation of the role of $\sigma$ and $\theta$).
For convenience, we denote $H_{j}(t)$ by $H(t)$.

We need to evaluate values $H(t)$ for random indices $t\in\{0,\ldots,N-1\}$.
Note that the signal $H$ results from the convolutions of two finite
bandwidth Box-car filters with the original signal $S$. Therefore,
any missing point needed by the two convolutions would lead to a failure
of computing $F(t)$. The total number of signal points involved depends
on the number of nonzero taps in these two filters. Moreover, random
dilation and modulation factors of the second Box-car filter make
computation more tricky.

One naive way is to dive into the two convolutions and sample each
signal point. If it is not available, stop evaluating this $F(t)$
and start with a new index $t$. This definitely increases time cost
by wasting abundant computation. For example, suppose five data are
needed and only one of them is missing, then the algorithm may compute
four data in vain in the worst case, where the missing data point
is visited last in the sequence of 5.

To avoid the above situation, we first compute the locations of all
the points that will be needed for the convolution; only if they are
all available will we start the computation. The locations related
to the convolution are given in the following lemma. 

\begin{lemma}  \label{lm:location}
Suppose we have a signal \(H(t)=( \chi_1^{(\sigma_1)} * ( \chi_{q_1}^{(\sigma_2)} * S)^{(\sigma_3)} )^{(\sigma_4)})(t)\), where \(\sigma_1\), \(\sigma_2\), \(\sigma_3\), and \(\sigma_4\) are dilation factors. From the definition of Box car filter, the taps for \(\chi_1\) lies in the interval \([-1, 1]\), the taps for \(\chi_{q_1}\) in \([-q_1, q_1]\), then in order to evaluate \(H(t)\), we need values of \(S\) with indices at \(\sigma_3 \sigma_4 t - \sigma_3 \sigma_1 i - j \sigma_2\), where integers \(i=-1,\ldots, 1\), \(j=-q_1,\ldots,q_1\).
\begin{proof}
To evaluate H(t), first let signal \( r=( \chi_{q_1}^{(\sigma_2)} * S)^{(\sigma_3)} \),
  then
    \begin{equation}
    H(t)=(\chi_1^{(\sigma_1)} * r)^{(\sigma_4)}(t) =\sum_{i=-1}^{ 1} \chi_1(\sigma_1 i) r(\sigma_4 t-\sigma_1 i)
    \end{equation}
     \begin{eqnarray}
    r(\sigma_4 t-\sigma_1 i)=( \chi_{q_1}^{(\sigma_2)} * S)^{(\sigma_3)}(\sigma_4 t-\sigma_1 i)
    = ( \chi_{q_1}^{(\sigma_2)} * S)(\sigma_3 \sigma_4 t -  \sigma_3 \sigma_1 i) \nonumber \\
        = \sum_{j=-q_1}^{q_1} \chi_{q_1}(\sigma_2 j) S(\sigma_3 \sigma_4 t- \sigma_3 \sigma_1 i - \sigma_2 j).
    \end{eqnarray}
   Thus, in order to get the value of \(H(t)\), we need values of all \(S(t^{'})\), where \(t^{'}=\sigma_3 \sigma_4 t- \sigma_3 \sigma_1 i - \sigma_2 j\), with \(i=-1,\ldots, 1\)
   and \(j=-q_1,\ldots, q_1\).
      \qquad\end{proof}
\end{lemma}\addtocounter{algorithm}{+1}

The scheme of the norm estimation algorithm is as follows.

\label{alg:norm}\begin{algorithm} {\textsc{Norm Estimation}} \label{estnorm} \\
Input: signal \(H\), \(k=0\), the number of iterations \(M=1.2\ln(1/\delta)\).\\
While \(k<M\):
\begin{enumerate}
    \item  Randomly generate the index \(t_k\).
    \item  Compute all indices needed by the two convolutions:
    \(\Upsilon=\{t^{'}, t^{'}=\sigma_3 \sigma_4 t- \sigma_3 \sigma_1 i - \sigma_2 j \}\), where \(i=-1,\ldots, 1\)
   and \(j=-q_1,\ldots, q_1\).
   \item   If all the points \( t^{'}\in \Upsilon\) are available, then compute \(H(t_k)\)
     else go to step 1 and generate another index \(t_k\).
    \item   estimate = 60-th percentile of the sequence \(\{|H(t_k)|^2 N\} \), where \(k=0,\ldots,M-1\).
   \end{enumerate}
\label{estnorm}\end{algorithm}\addtocounter{theorem}{+1}

If there exist satisfactory data groups, although maybe very few,
the norm estimation will eventually find them. However, when most
data are unavailable, the program may struggle in a long loop and
take a huge amount of time. We introduce some tricks to avoid this.
For example, set an upper bound MAX on the number of the loops. If
it is reached, just use the sample points generated so far to estimate
the norms. This technique may lead to a larger error, and thus hamper
our frequency identification. However, by repeating the calculation,
as stipulated by Lemma 3.2, we reduce the inaccuracy. Anyway we cannot
hope to recover the signal, if $p$ is too small.

The following lemma investigates the number of repetitions to get
a satisfactory data group for estimating norms.

\begin{lemma} 
Suppose \(\chi_{q_1}\) and \(\chi_{q_2}\) are two Box-car filters with numbers of taps \(2q_1+1\) and \(2q_2+1\) respectively. Define \(D_{q_1,q_2} = \chi_{q_1} * \chi_{q_2}\). Then $D_{q_1,q_2}$ has $2q_1+2q_2+1$ nonzero taps in the time domain.
\end{lemma}\addtocounter{algorithm}{+1}

\begin{lemma}
Randomly choose an index for signal \(H(t)\), then after \(k>\log \delta /
\log(1-(1-p)^{2q_1+2q_2+1})\) iterations, we can get at least one satisfactory index with high probability
\(1-\delta\).
\end{lemma}\addtocounter{algorithm}{+1} \begin{proof}
It is easy to prove by Lemma \ref{lm:rept}.
 \qquad\end{proof}

Here is a new scheme for estimating norms, which uses much fewer samples
than the original one and still achieves good estimation. In \cite{Zou},
we propose a lemma that enabled us to achieve a good norm estimation
by only a few samples. The following lemma is its adaption to the
case of unevenly spaced data.

\begin{lemma}
If a signal \(H\) is 95\% pure and if \(r>1.2 \ln (1/\delta)\), the output of Algorithm \ref{estnorm} gives an estimation of its energy
which exceeds \(\|H\|^2/3\) with probability exceeding \(1-\delta\).
\end{lemma}\addtocounter{algorithm}{+1} \begin{proof}
The proof is very similar to that of Lemma 4.5 in \cite{Zou}. We shall present only the difference of these two proofs. Suppose we sample \(r\) times for the signal \(H\). Let \(\kappa=\{t:N|H(t)|^2<\|H\|^2/3 \}\), with \(\kappa^c\) as its complement, we have
\begin{equation}
\left |\sum_{t \in \kappa}H(t) \right |^2 \leq |\kappa| \sum_{t \in \kappa}|H(t)|^2 \leq |\kappa|^2 \frac{1}{N}\frac{1}{3} \|H\|^2.
\end{equation}
 On the other hand, we know that the signal is 95\(\%\) pure, i.e. \(|\hat{H}(\omega_0)|^2 \geq 0.95\|H\|^2\) for some \(\omega_0\). By modulating, \(\omega_0\) can be moved to 0; therefore, we can, without loss of
generality, suppose most of the energy concentrates at the frequency 0; then
\begin{equation}
\left |\frac{1}{\sqrt{N}} \sum_{t=1}^N H(t) \right |^2 = |\hat{H}(0)|^2 \geq 0.95 \|H\|^2.
\end{equation}
So we have 
\begin{eqnarray}
\left |\sum_{t \in \kappa^C}H(t) \right | \geq \sqrt{0.95N} \|H\| - |\kappa| \frac{1}{\sqrt{3N} } \|H\|.
\end{eqnarray}
On the other hand,\(|\sum_{t\in\kappa^{C}}H(t)|\leq|\kappa^{C}|\| H\|=(N-|\kappa|)\| H\|\), so that
\begin{equation}
N-|\kappa| \geq \left (\sqrt{0.95N} - \frac{|\kappa|}{\sqrt{3N}} \right )^2.
\end{equation}
Let \(\alpha = \frac{|\kappa|}{N}\); the above inequality becomes
\begin{equation}
\alpha^2 + \left( 3-2 \sqrt{0.95*3}\right) \alpha -0.15 \leq 0.
\end{equation}
Thus \(0 \leq \alpha \leq 0.075 \). 
Define now a random variable \(X_{\kappa}= \left (\sum_{i=1}^N \chi_{\kappa}(i) \right )\); it will be useful to estimate
\begin{equation}
E(X_{\kappa})=\frac{|\kappa|}{N} \leq 0.075,
\end{equation}
and the expectation of the random variable \(e^{z X_{\kappa}}\),
\begin{equation}
E(e^{X_{\kappa} z}) = e^0 Prob(\chi_{\kappa}(i)=0) + e^z Prob(\chi_{\kappa}(i)=1) = 1-\alpha + \alpha e^z.
\end{equation}
Suppose now we sample the signal \(H\) \(r\) times, and take the \mbox{\textit{{60-th}}} percentile of the numbers \(N|H(t_1)|^2, \ldots, N|H(t_r)|^2\).
By Chernoff's standard argument and similar procedure of Lemma 4.5 in \cite{Zou}, we have for \(z>0\),
\begin{eqnarray*}
Prob \left (\mbox{\textit{{60-th}}} \, percentile < \frac{1}{3} \|H\|^2 \right ) = \left [ (1-\alpha) e^{-0.6z} + \alpha e^{0.4 z} \right ]^r. \nonumber 
\end{eqnarray*}
Take \(z=\ln (1.25(1-\alpha)/ \alpha)\), then
\begin{equation}
(1-\alpha)e^{-0.6z} + \alpha e^{0.4z} = 1.97 \alpha^{0.6} (1-\alpha)^{0.4}.
\end{equation}
The right hand side of (35) is increasing in \(\alpha\) on the interval \([0, 0.075]\); since \(\alpha \leq 0.075\), we obtain an upper bound by substituting \(0.075\) for \(\alpha\):
\begin{eqnarray}
\left [ (1-\alpha) e^{-0.6z} + \alpha e^{0.4 z} \right ]^r = \left [ 1.97 \alpha^{0.6} (1-\alpha)^{0.4} \right ]^r \leq e^{-0.90 r}.
\end{eqnarray}

For \(Prob \left (\mbox{\textit{{60-th}}} \, percentile < \frac{1}{3} \|H\|^2 \right ) \leq \delta\), we need \(r \geq 1.2 \ln (1/\delta)\), we have
\begin{equation}
Prob(Output \geq \|H\|^2/3) = Prob(\mbox{\textit{{60-th}}}\, percentile\,of\, N|H(t)|^2 \geq \|H\|^2/3) \geq 1-\delta.
\end{equation}
 \qquad\end{proof}

This norm estimation procedure will be used repeatedly in the group
testing step below.

\subsection{Isolation}

For a significant frequency in signal $S$, isolation aims to construct
a series of new signals, such that this significant frequency becomes
predominant in at least one of the new isolation signals.

\begin{lemma}
Given signals \(S\), \(S_1\), and the parameters as stated in Lemma \ref{lm:smallfilter}. Suppose \(F_1(t) = S_1(t)*\chi_1(t) = (\chi_T(t) S(t))*\chi_1(t)\), \(F(t)= S(t) * \chi_1(t)\). If \( p\geq \frac{1}{1+(N-1)\lambda \tau^2}\), then for each \(\omega\) with \(|\hat{S}(\omega)|^2 > B \lambda \|S\|^2\), isolation algorithm can create a signal $F_1^{*}$, such that
 \begin{equation}
|\hat{F}_1^{*}(\omega)|^2 \geq 0.98\|F_1^{*}\|^2.
\end{equation}
 \end{lemma}\addtocounter{algorithm}{+1}\label{lem:iso} \begin{proof}
Since \( |\hat{S}(\omega)|^2 > B \lambda \|S\|^2\), we have \(|\hat{S}(\omega)| > \sqrt{B \lambda} \|S\|\). Then there exists some \(\eta>0\), such that \(|\hat{S}(\omega)| \geq (\sqrt{\eta}+\sqrt{B \lambda})\|S\|.\) 
 Lemma \ref{lm:diffest} states that 
\(|\hat{S}_1(\omega)-\hat{S}(\omega)| \leq \sqrt{B \lambda} \|S\|\). Therefore
\begin{equation}
|\hat{S}_1(\omega)|\geq \sqrt{\eta}\|S\|\geq \sqrt{\eta}\|S_1\|.
\end{equation}
Isolation algorithm returns \(F_1^{(0)}, \ldots, F_1^{(2k)}\) with \(k<O(\frac{1}{\eta})\), as described in \cite{GGIMS}. For any \(\omega\) with \(|\hat{S}_1(\omega)|^2 \geq \eta \|S_1\|^2\), there exists some \(j\), such that
\begin{equation}
|\hat{F}_1^{(j)}(\omega)|^2 \geq 0.98 \|F_1^{(j)}\|^2.
\end{equation}

Let $F_1^{*}=F_1^{(j)}$, then  

\begin{equation}
|\hat{F}_1^{*}(\omega)|^2 \geq 0.98\|F_1^{*}\|^2.
\end{equation}

\end{proof}

Theoretically, in order to capture a significant mode, we need $O(1/\eta)$
signals. However, in practice, much fewer signals is enough to achieve
this goal.

\subsection{Group Testing}

Isolation has produced several signals, one of which contains the
most significant frequency. Group testing uses repeated zoom-ins on
one of the signals, and norm testing to select where to zoom in, in
order to determine the frequency. The goal of group testing is thus
to find the most significant mode of the signal $F_{1}^{*}$ from
isolation. It uses recursive procedures MSB (Most Significant Bit)
to approach this mode gradually.

\emph{Definition}: Denote a set $\{\omega:\,\,(2l-1)N/32\leq\omega\leq(2l+1)N/32\}$
by $interval_{l}$.

Group test algorithm is given as follows.

\begin{algorithm} \label{alg:grouptest}{\textsc{Group Testing}} \\
Input isolation signal \(F_1^{*}\) to \(F_1^{(0)}\), \(i=0\), \(q=1\) \\
While \(q<N\), in the \(i\)-th iteration,
 \begin{enumerate}
       \item  Find the most significant bit \(v\) and the number of significant intervals \(c\) by the procedure MSB.
    \item  Update \(i=i+1\), modulate the signal \(F_1^{(i)}\) by \( \lfloor (v+0.5)N/16 \rfloor \) and dilate it by a factor of \( \lfloor 16/c \rfloor\). Store it in \(F_1^{(i+1)}\).
    \item  Call Group Test again with the new signal \(F_1^{(i)}\), denote its output by \(g\).
    \item  Update the accumulation factor \(q = q * \lfloor 16/c \rfloor \).
      \item If \(g> N/2\), then \(g = g -N\).
      \item return \( \lfloor g/\lfloor 16/c \rfloor + (v+1/2)N/16+0.5 \rfloor(mod\,\, N)\);
     \end{enumerate}
 \end{algorithm}\addtocounter{theorem}{+1}

The MSB procedure is as follows.

\begin{algorithm} \label{alg:msb}{\textsc{MSB (Most Significant Bit)}} \\
\text{}\hspace{10mm} Input: signal \(F_1^{(i)}\) with length \(N\), a threshold \(0<\eta<1\).
\begin{enumerate}
    \item Get a series of new signals \(H_j(t) =F_1^{(i)}(t) \star (e^{2 \pi i j t/16} \chi_1 )\), \(j=0, \ldots, 15\). 
    \item Estimate the energies \(e_j\) of \(H_j\), \(j=0, \ldots, 15\).
    \item for \(l=0,\ldots,15\), compare the energies \(e_l\) with all other energies \(e_j\), where \(j=(l+4)mod\,16, (l+5)mod\,16, \ldots,(l+12)mod\,16\). If \(e_l > e_j\) for all these \(j\), label it as an interval with large energy.
     \item Find the longest consecutive intervals of large energies. Take their center as \(v \), and the number of those intervals as \(c \).
    \item If \(c<8\), then do the original MSB in {\emph{\cite{GGIMS}}} to get \(v\) and set \(c=8\);
    \item Return the dilation-related factor \(c\) and the most significant bit \(v\).
\end{enumerate} \label{alg:msb}
\end{algorithm}\addtocounter{theorem}{+1}

For convenience, we denote $F_{1}^{(i)}$ by $\mathbf{F_{1}}$.

\begin{lemma}
Given a \( 98\%\) pure signal \(\mathbf{F_1}\), suppose \(G_j(t) = e^{2 \pi i j t /16} \chi_1(t)\). Then Algorithm \ref{alg:grouptest}, with Algorithm \ref{alg:msb} as its subroutine, can find the significant frequency \(\omega_1\) of the signal \(\mathbf{F_1}\) with high probability.
\end{lemma}\addtocounter{algorithm}{+1} \begin{proof}
The proof is similar to that of Lemma 5 in \cite{GGIMS}, with some changes: 

Since the signal \(\mathbf{F_1}\) is \(98\%\) pure, there exist a frequency mode \(\omega_1\) and a signal \(\rho\), such that \(\mathbf{F_1}=a\phi_{\omega_1}+\rho\), where $|a|^2 \geq 0.98\|\mathbf{F_1}\|^2$ and \(\|\rho\|^2 \leq 0.02\|\mathbf{F_1}\|^2\). Without loss of generality, assume \(\omega_1 \in [-N/32, N/32]\). The whole region is divided into 16 subintervals \([jN/16-N/32, jN/16+N/32]\), where \(j=0,\ldots, 15\). To estimate \(\widehat{\mathbf{F_1}*G_0}(\omega_1)\) for \(|\omega_1|\leq N/32\), we use that \(|\hat{G}_0(\omega_1)|=|\hat{\chi}_1(\omega_1)|\geq 0.987\) for \(|\omega_1|\leq N/32\). It follows that 
\begin{eqnarray*}
|\widehat{\mathbf{F_1} \ast G_{0}}(\omega_1)|^{2}  = N \left|\hat{\mathbf{F}}_1(\omega_1)\hat{G}_0(\omega_1)\right|^2 \geq  N 0.987^{2}|\hat{\mathbf{F}}_1(\omega_1)|^{2} \geq N 0.987^{2}0.98\| \mathbf{F_1}\|^{2} \nonumber \\
\geq 0.954N\|\hat{\mathbf{F}}_1\|^{2}  \geq 0.954N\|\hat{\mathbf{F}}_1\hat{G_{0}}\|^{2}=0.954\|\mathbf{F_1}\ast G_{0}\|^{2}.
\end{eqnarray*}
Therefore the estimation \(X\) of \( \|\mathbf{F_1} * G_0\|\) satisfies:
 \begin{eqnarray*}
X \geq \|\mathbf{F_1} * G_0\|^2/3 = \|\widehat{\mathbf{F_1}*G_0}\|^2/3  = \sum_{\omega} |\widehat{\mathbf{F_1}*G_0}(\omega)|^2 /3 \geq |\widehat{\mathbf{F_1}*G_0}(\omega_1)|^2/3 \nonumber \\
\geq 0.954N\|\mathbf{F_1}\|^2/3 \geq 0.318 N\|\mathbf{F_1}\|^2.
\end{eqnarray*}
Next consider the energy of \(\mathbf{F_1}*G_{4}\).
 \begin{eqnarray*}
\|\hat{\rho}\hat{G_{4}}\|^{2} =\sum_{\omega}|\hat{\rho}(\omega)\hat{G_{4}}(\omega)|^{2} \nonumber \\
\leq \sum_{\omega}|\hat{\rho}(\omega)|^{2} =  \| \rho\|^{2}\leq 0.02 \| \mathbf{F_1}\|^{2}.
\end{eqnarray*}
 Since $|\hat{G}_4(\omega_1)|<0.464$, we have
\begin{eqnarray*}
|\hat{\mathbf{F}}_1(\omega_1)\hat{G}_{4}(\ \omega_1)|\leq|\hat{\mathbf{F}}_1(\omega_1)||\hat{G}_{4}(\ \omega_1)| \leq |\hat{\mathbf{F}}_1(\omega_1)|0.464 \leq 0.464 \| \mathbf{F_1}\|
\end{eqnarray*}
Also \( \|\hat{\mathbf{F}}_1 \hat{G}_{4}\|^{2}-|\hat{\mathbf{F}}_1(\omega_1) \hat{G}_{4}(\omega_1)|^{2} \leq 0.02\| \mathbf{F_1}\|^{2}\).
Thus

\[ \| \hat{\mathbf{F}}_1 \hat{G}_4 \|^2 \leq 0.464^2\| \mathbf{F_1}\|^{2}+0.02\| \mathbf{F_1}\|^{2}=0.24\| \mathbf{F_1}\|^{2}. \] It follows that \[ \| \mathbf{F_1} \ast G_4 \|^2 =\| \widehat{\mathbf{F_1}*G_4}\|^2 = N \|\hat{\mathbf{F}}_1 \hat{G}_4\| \leq 0.24N\| \mathbf{F_1}\|^{2}. \]
Then we compare \(\| \mathbf{F_1}\ast G_{4}\|^{2}\) with the lower bound of the estimation of \(\| \mathbf{F_1}\ast G_{0}\|^{2}\), which is \[ 0.24N \| \mathbf{F_1}\|^{2} \leq 0.318N \| \mathbf{F_1}\|^{2}, \]  which is less than the estimation for \(\| \mathbf{F_1}\ast G_{0}\|^{2}.\) In general, \(\omega\in interval_{j}\), for \(j\) not necessarily 0. Therefore we compare \(\| \mathbf{F_1}\ast
G_{j^{'}}\|^{2}\)with \(\| \mathbf{F_1}\ast G_{j}\|^{2}\), where \(|j-j^{'}|\geq4\). If there is some \(j\) with \(\| \mathbf{F_1}\ast G_{j}\|^{2}\) apparently
larger than \(\| \mathbf{F_1}\ast G_{j^{'}}\|^{2}\), then we conclude \(\omega_1 \notin interval_{j^{'}}\). Otherwise, possibly
\(\omega_1\in interval_{j^{'}}\). By the above argument, we can always eliminate 9 consecutive interval regions out of 16, leaving a cyclic interval of length at most \(7N/16\). The remaining proof is exactly the same as Lemma 8 in paper \cite{GGIMS}.
\end{proof}

Remark: In \cite{Zou}, we showed that group testing works for a Box-car
filter with width more than $21$, i.e. $k>10$. In that case, $2k+1$
intervals are sufficient. A similar conclusion still holds in the
unevenly spaced data case. However, the lemma above proves the success
of group testing under different conditions. In our proof, we use
a Box-car filter with much shorter width, namely 3 in time domain;
this works well if 16 intervals are taken. In practice, we use these
shorter filters; we can usually (if $B$ is small) get away with using
much fewer intervals as well (e.g. 3 instead of 16).

\subsection{Adaptive Greedy Pursuit}

In summary, given a signal $S$, for an accuracy $\epsilon$ and for
$B$ modes, we can find a very good approximation of the signal $S$
by using Algorithm \ref{alg:total1}. 

\begin{theorem}
\label{lm:totalcost} Given a signal \(S\), an accuracy \(\epsilon\), success probability \(1-\delta\), Algorithm \ref{alg:total1} can output a \(B\)-term representation \(R\) with sum-square-error \(\|S-R\|^2\leq (1+\epsilon) \|S-R_{opt}\|^2\), where \(R_{opt}\) is the \(B\)-term representation for \(S\) with the least sum-square-error, with time and space cost \(poly(B,\log(N), \frac{1}{\epsilon}, \log(1/\delta))\) for computing and \(\frac{B \log M \log N \log \delta}{\lambda log(1-(1-p)^{2q_1+2q_2+1})}\) \(+\frac{\log (1/\delta) \log M}{\lambda \log p}\) for just visiting samples. 
\begin{proof}
{\textrm{We omit the proof since it is very similar to Theorem 9 in \cite{GGIMS}.}}
\end{proof}
\end{theorem}\addtocounter{algorithm}{+1}

\section{NERA$\ell$SFA with Interpolation Technique}

The greedy algorithm described above is fast. When $p$ is sufficiently
large (e.g. $p>0.7$), the approach proposed and discussed in the
previous section works well. For smaller $p$, the amount of time
wasted to find available sample groups becomes unacceptably long.
For example, when $B=2$, $N=100$, $p=0.4$, the algorithm couldn't
find the signal within 200 greedy pursuit iterations. For this reason,
we introduced an interpolation technique to get an approximate value
of the missing point in the norm estimation procedure. This algorithm
is efficient even in smaller $p$ cases.

\subsection{Lagrange Interpolation Technique}

The task of interpolation is to estimate $S(t)$ for arbitrary $t$
by drawing a smooth curve through all the known points \cite{PTVF}.
It is called interpolation when the desired $t$ is between the largest
and smallest of these $t_{i}$'s. We use Lagrange Polynomial Interpolation,
one of the simplest and most popular interpolation techniques.

Generally, the number of interpolation points determines the degree
of a polynomial. A polynomial of higher degree is smoother with smaller
approximation errors at the expense of more computation. Thus we choose
a second degree polynomial, as a balance between computational complexity
and accuracy. It is given explicitly by Lagrange's classical formula.
If the three nearest neighbors are $(t_{1},S(t_{1}))$, $(t_{2},S(t_{2}))$,
$(t_{3},S(t_{3}))$, the polynomial is \begin{equation}
P(t)=\frac{(t-t_{2})(t-t_{3})}{(t_{1}-t_{2})(t_{1}-t_{3})}S(t_{1})+\frac{(t-t_{1})(t-t_{3})}{(t_{2}-t_{1})(t_{2}-t_{3})}S(t_{2})+\frac{(t-t_{2})(t-t_{1})}{(t_{3}-t_{2})(t_{3}-t_{1})}S(t_{3})\end{equation}

If $S(t)$ is three times differentiable in an interval $[a,b]$,
and the points $t_{1},t_{2},t_{3}\in[a,b]$ are different, then there
exists some $v\in[a,b]$, such that the approximation error is $S(t)-P(t)=\frac{S^{(3)}(v)}{3!}(t-t_{1})(t-t_{2})(t-t_{3})$.

\subsection{Estimate Norms with Interpolation}

We introduce the interpolation scheme into estimating norms. The idea
is to estimate the value of a missing point by the Lagrange interpolation.
The detailed algorithm for estimating norms is as follows.

\begin{algorithm}  {\textsc{Estimate Norm with interpolation technique}} \\
Input: signal \(H\), \(k=0\), the maximum number of samples \(M\).
 \begin{enumerate}
    \item  Randomly generate the index \(t_k\), where \(k=0,\ldots, M-1\).
      \item  For each $k$, if \(H(t_k)\) is not available, estimate \( H(t_k) \) by Lagrange interpolation; else compute \(H(t_k)\) directly.
    \item   Estimation = 60-th percentile of the sequence \(\{|H(t_k)|^2 N\} \), where \(k=0,\ldots,M-1\).
   \end{enumerate}
\end{algorithm}

Note that we use interpolation \emph{only} in norm estimation steps,
where precision is less critical. With less precise norm estimation,
the localization of important modes could still work well when iterated.
For coefficient estimation, which needs to be more precise, we always
search for available samples.

\section{Numerical Results}

In this section, we present striking numerical results of NERA$\ell$SFA,
comparing to the Inverse Non-equispaced Fast Fourier Transform (INFFT)
algorithms. The popular benchmark software NFFT version 2.0 is used
to give performance of INFFT, with default CGNE\_R method and Dirichlet
kernel. Its time cost excludes the precomputation of samples values,
which takes $O(L)$. Numerical experiments show the advantage of our
NERA$\ell$SFA algorithm in processing large amount of data. We begin
in Section 5.1 with comparing NERA$\ell$SFA with INFFT for some one
and two dimensional examples with different length. In Section 5.2,
the performance for different number of modes is shown. Finally, we
test the capability of NERA$\ell$SFA to recover the signal in the
situation with a large amount of missing data and in presence of large
noise.

All the experiments were run on an AMD Athlon(TM) XP1900+ machine
with Cache size 256KB, total memory 512 MB, Linux kernel version 2.4.20-20.9
and compiler gcc version 3.2.2. The numerical data is an average of
10 runs of the code; errors are given in the $L^{2}$ norm.

\subsection{Experiments with Different Length of Signals}

We ran the comparison for a 8-mode superposition signal $S(t)=\sum_{i=1}^{B}\phi_{\omega_{i}}$,
plus white noise $\nu$ with the standard deviation $\sigma=0.5$,
damped by a factor of $1/\sqrt{N}$, ( so that $\Vert\nu\Vert^{2}=\sigma^{2}=0.25$;
since $\Vert S\Vert^{2}=8$, this implies $SNR=20\log_{10}32\thickapprox30.1dB$).
Other parameters are $B=8$, $\epsilon=0.02$, $\delta=0.01$, and
$p=70\%$. The missing data are randomly and uniformly distributed.
NERA$\ell$SFA outperforms INFFT in speed when $N$ is large; see
Table \ref{tab:B13} and Figure \ref{fig:diffN1d}. The corresponding
crossover point is $N\geq2^{15}=32768$ . For example, to process
$2^{19}=524,288$ data, more than nineteen minutes (estimated) are
needed for INFFT versus approximately one second for NERA$\ell$SFA.
Experiments support the theoretical conclusion that NERA$\ell$SFA
would be faster than INFFT after some $N$ for a sparse signal; whatever
the sparsity, i.e. whatever the value of $B$, there always exists
some crossover $N$.

\begin{table}
\begin{center}\begin{tabular}{|c|c|c|c|}
\hline 
N &
 INFFT&
 NERA$\ell$SFA &
 NERA$\ell$SFA \tabularnewline
&
&
 (+sampling)&
 (w/o sampling)\tabularnewline
\hline
$2^{9}$=512&
 0.01&
 0.63&
 0.31\tabularnewline
$2^{11}$=2048&
 0.03&
 0.77&
 0.37\tabularnewline
$2^{13}$=8192&
 0.17&
 0.90 &
 0.46\tabularnewline
$2^{15}$=32768&
 0.83&
 0.93 &
 0.49\tabularnewline
$2^{17}$=131072&
 4.30&
 1.03&
 0.51\tabularnewline
$2^{19}$=524288&
 19.94&
 1.20&
 0.61  \tabularnewline
\hline
\end{tabular}\end{center}

\caption{\label{tab:B13}Experiments with fixed $B=8$, $p=0.7$, $d=1$ (one
dimension), and varying length $N$ of signals; an i.i.d. white noise
is added with $\sigma=0.5$, or $SNR\simeq30dB$ (see text). For each
length of the signal, 10 different runs were carried out; the average
result is shown. We did all the tests for NERA$\ell$SFA with Lagrange
interpolation, as explained in the text. Two kinds of time costs for
NERA$\ell$SFA are provided. One is the total running time and another
is the running time excluding the sampling time. The time of INFFT
does not include the precomputation time for samples. }
\end{table}

\addtocounter{figure}{+1}

\begin{figure}
\begin{center}\includegraphics[%
  width=10cm]{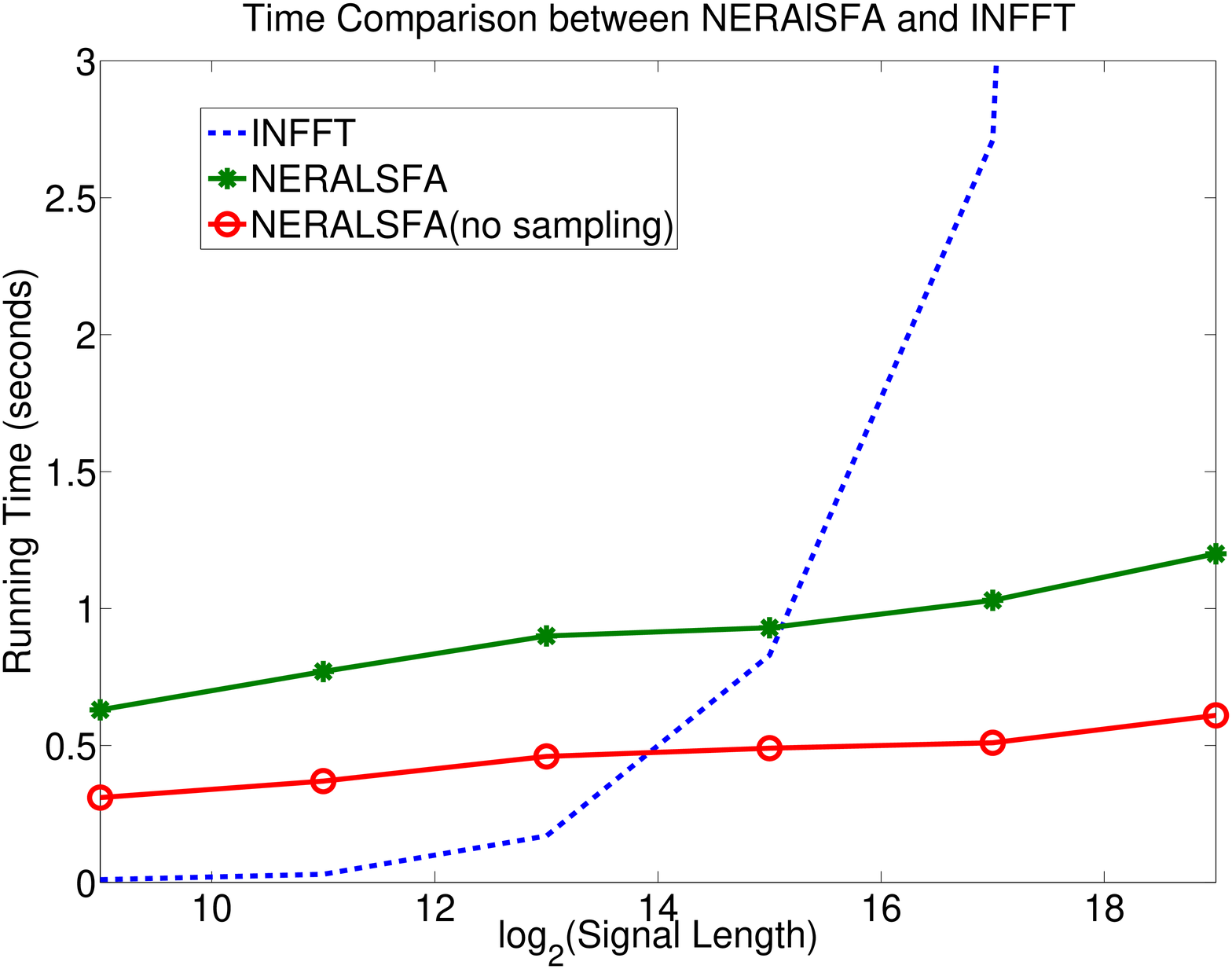}\end{center}

\caption{Time Comparison between INFFT and NERA$\ell$SFA for different $N$
with $B=8$, $p=0.7$, $d=1$. The result in Table \ref{tab:B13}
is shown in the form of a graph here. The $x$ coordinate is the $\log_{2}(N)$,
the $y$ coordinate presents the running time for each algorithm.
NERA$\ell$SFA without sampling surpasses INFFT at $N=2^{14}=16384$.
\label{fig:diffN1d}}
\end{figure}

\addtocounter{table}{+1}

In two dimensions, we test a noisy 6-mode superposition signal $S(t)=\sum_{i=1}^{B}\phi_{\omega_{xi}}\phi_{\omega_{yi}}+\nu$,
with $B=6$, $\epsilon=0.02$, $\delta=0.01$, $p=80\%$, and $\sigma=0.1$.
Missing data are randomly and uniformly distributed. As the number
of grid points $N$ in each dimension grows, two dimensional NERA$\ell$SFA
outperforms two dimensional INFFT at $N\geq512$, as Table \ref{tab:2dB13}
and Figure \ref{fig:diffN2d} show. The crossover point becomes much
smaller in high dimensions situation. It would not be surprising that
for recovering a 6-mode three dimensional signal, NERA$\ell$SFA surpasses
INFFT at a hundred sampling grid points in each dimension.

\begin{table}
\begin{center}\begin{tabular}{|c|c|c|c|}
\hline 
N &
 INFFT&
 NERA$\ell$SFA &
 NERA$\ell$SFA \tabularnewline
&
&
 (+sampling)&
 (w/o sampling)\tabularnewline
\hline
$128$&
 0.13&
 2.86 &
 1.57\tabularnewline
256&
 0.73&
 2.60&
 1.46\tabularnewline
512&
 3.00&
 3.70&
 2.13\tabularnewline
1024&
 11.59&
 4.31&
 2.94\tabularnewline
$2048$&
 54.94&
 6.56&
 4.90  \tabularnewline
\hline
\end{tabular}\end{center}

\caption{\label{tab:2dB13}Experiments with fixed $B=6$, $p=0.8$, $d=2$
(two dimensions), and varying length $N$ of signals; an i.i.d white
noise is added with $\sigma=0.1$, or $SNR\simeq56dB$ (see text).
For each length of the signal, 10 different runs were carried out;
the average result is shown. We did all the tests for NERA$\ell$SFA
with two dimensional interpolation techniques as shown in the appendix.
Again, two kinds of time costs for NERA$\ell$SFA, the one with and
without sampling time is provided. The time of INFFT excludes the
sampling time. }
\end{table}

\addtocounter{figure}{+1}

\begin{figure}
\begin{center}\includegraphics[%
  width=10cm]{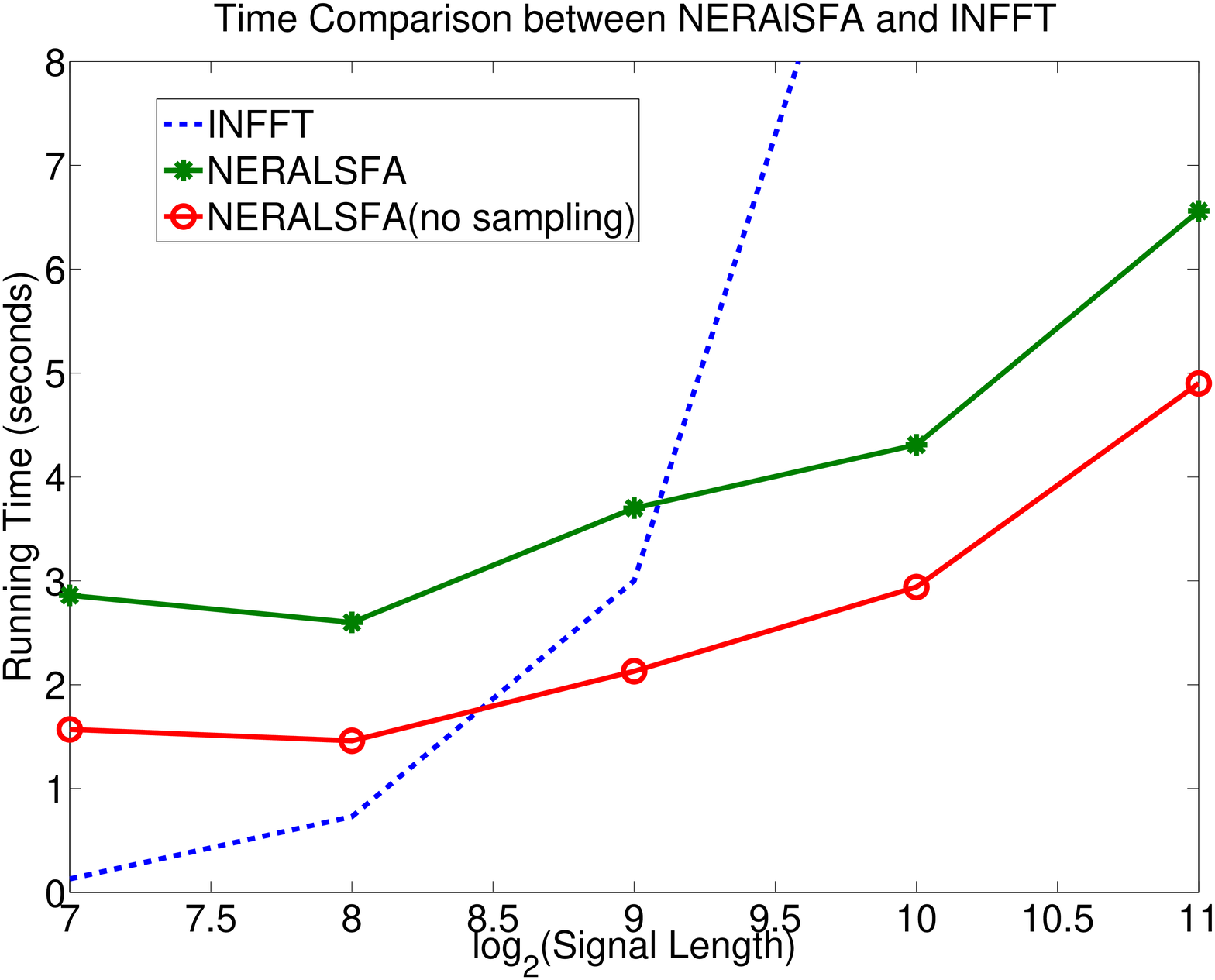}\end{center}

\caption{\label{fig:diffN2d}Time comparison between INFFT and NERA$\ell$SFA
for different $N$ with fixed $B=6$, $p=0.8$, $d=2$. The $x$ coordinate
is the logarithm of length $N$ of signal in each dimension. INFFT
is very fast when $N$ is relatively small and slows down quickly
as $N$ increases. On the contrary, it takes NERA$\ell$SFA similar
time to process small and large $N$ problem. NERA$\ell$SFA without
sampling outperforms INFFT at $N=2^{8.5}$=362. }
\end{figure}

\addtocounter{table}{+1}

\subsection{Experiments with Different Number of Modes}

The number of modes has an important influence on the running time
since the crossover point varies for signals with different $B$.
To investigate this, we did the experiments with fixed $N=2^{18}=262144$,
$p=0.6$ and varying $B$. As before, we take $S$ to be a superposition
of exactly $B$ modes with white noise, i.e. $S(t)=\sum_{i=1}^{B}c_{i}\phi_{\omega_{i}}+\nu$,
with standard deviation of noise $\sigma=0.05$. Available data are
uniformly and randomly distributed. Table \ref{tab:diffB} and Figure
\ref{fig:diffB1d} compare the running time for different $B$ using
INFFT and NERA$\ell$SFA. At first, NERA$\ell$SFA takes less time
because $N$ is so large. However, the execution time of INFFT keeps
constant for different number of modes $B$, while that of modified
RA$\ell$SFA is polynomial of higher order. INFFT is faster than NERA$\ell$SFA
when $B\geq10$. The regression techniques shows empirically that
the order of $B$ in NERA$\ell$SFA is greater than quadratic. This
is one of the characteristics of this version of the RA$\ell$SFA
algorithms and irrelevant to the nonequispaceness of the data. (A
different version of RA$\ell$SFA in \cite{GMS} is linear in $B$,
but maybe less easily used when not all equispaced data are available.
)

\begin{table}
\begin{center}\begin{tabular}{|c|c|c|c|c|}
\hline 
number of modes&
 SNR&
 NERA$\ell$SFA &
 NERA$\ell$SFA &
 INFFT \tabularnewline
$B$&
 (dB)&
 (+sampling) &
 (w/o sampling)&
\tabularnewline
\hline
2&
 58&
 0.06&
 0.01&
 1.35\tabularnewline
\hline
4&
 64&
 0.24&
 0.06&
 1.35\tabularnewline
\hline
6&
 68&
 0.61&
 0.23&
 1.35\tabularnewline
\hline
8&
 70&
 1.44&
 0.69&
 1.35\tabularnewline
\hline
10&
 72&
 2.45&
 1.39&
 1.35\tabularnewline
\hline
13&
 74&
 5.78&
 3.64&
 1.35\tabularnewline
\hline
16&
 76&
 10.03&
 7.17&
 1.35  \tabularnewline
\hline
\end{tabular}\end{center}

\caption{\label{tab:diffB}Experiments with fixed $N=2^{18}$, $p=0.6$, $d=1$
(one dimension), $\sigma=0.05$, and varying number of modes $B$
of signals. For each length of the signal, 10 different runs were
carried out; the average result is shown. We did all the tests for
NERA$\ell$SFA with interpolation techniques. We present two different
time costs of NERA$\ell$SFA, with and without sampling. }
\end{table}

\addtocounter{figure}{+1}

\begin{figure}
\begin{center}\includegraphics[%
  width=10cm]{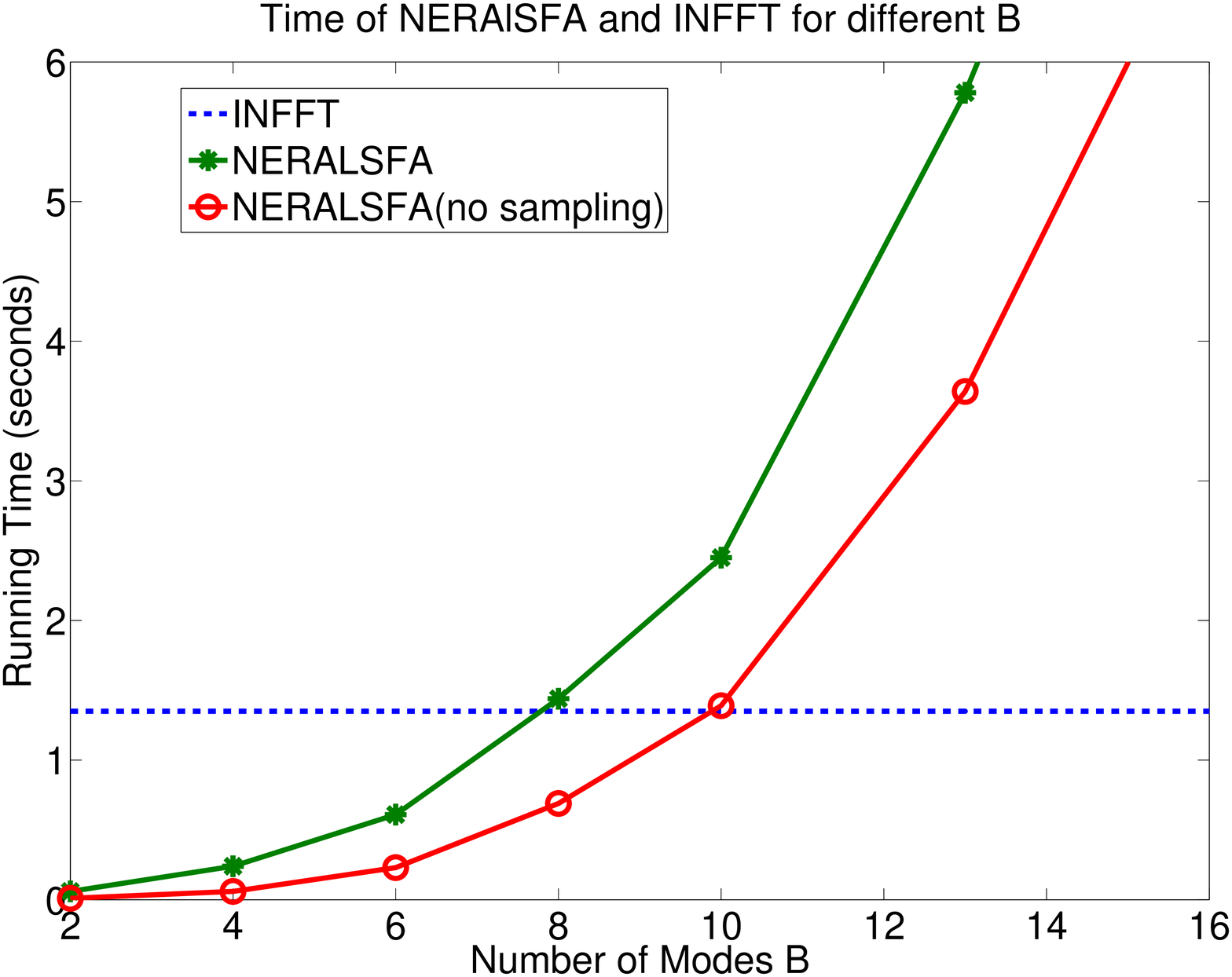}\end{center}

\caption{\label{fig:diffB1d}Time Comparison between INFFT and NERA$\ell$SFA
for different $B$ with with fixed $N=2^{18}$, $p=0.6$, $d=1$ (one
dimension), $\sigma=0.05$, a graph of the result in Table \ref{tab:diffB}.
The $x$ coordinate is the number of modes $B$, the $y$ coordinate
presents running time. The running time of NERA$\ell$SFA is polynomial
to $B$. In contrast, the time of INFFT keeps constant for different
$B$, excluding precomputation for the samples. NERA$\ell$SFA without
sampling begins to be slower than INFFT at $B=10$ for $N=2^{18}$.}
\end{figure}

\addtocounter{table}{+1}

\subsection{Experiments for Different Percentage of Missing Data}

The advantage of interpolation techniques is to recover a signal even
when a large percentage of data is missing. Table \ref{tab:B2rec}
shows the recovery effect for a two-mode pure signal $c_{1}\phi_{\omega_{1}}+c_{2}\phi_{\omega_{2}}$,
$N=10^{6}$ with all the other parameters $\epsilon$ and $\delta$
the same as before. When the percentage of available data is large,
both algorithms recover the signal well with similar running time.

\begin{table}
\begin{center}\begin{tabular}{|c|c|c|c|c|}
\hline 
p &
 Time of NERA$\ell$SFA&
 success&
 Time of NERA$\ell$SFA &
 success\tabularnewline
&
 (with interpolation)&
 probability&
 (w/o interpolation)&
 probability \tabularnewline
\hline
1 &
 0.03 &
 100 $\%$&
 0.03 &
 100 $\%$\tabularnewline
\hline
0.8 &
 0.04&
 100 $\%$&
 0.06 &
 100 $\%$\tabularnewline
\hline
0.6 &
 0.05 &
 100 $\%$&
 0.49 &
 100 $\%$\tabularnewline
\hline
0.4 &
 0.05&
 100 $\%$&
 0.45&
 100 $\%$\tabularnewline
\hline
0.3 &
 0.06 &
 100 $\%$&
 - &
 0 $\%$\tabularnewline
\hline
0.2 &
 0.06 &
 100 $\%$&
 - &
 0 $\%$\tabularnewline
\hline
0.1 &
 0.07 &
 100 $\%$&
 - &
 0 $\%$\tabularnewline
\hline
$10^{-2}$&
 0.11 &
 100 $\%$&
 - &
 0 $\%$\tabularnewline
\hline
$10^{-3}$&
 0.51 &
 100 $\%$&
 - &
 0 $\%$\tabularnewline
\hline
$10^{-4}$&
 4.58 &
 100 $\%$&
 - &
 0 $\%$\tabularnewline
\hline
$0.00002$&
 758.22 &
 97 $\%$&
 - &
 0 $\%$ \tabularnewline
\hline
\end{tabular}\end{center}

\caption{\label{tab:B2rec}Experiments with fixed $B=2$, $N=10^{6}$, no
noise, and varying percentage of available data. Each entry is based
on the average of 10 different runs. In each run, the number of iterations
is limited to 200; (this also corresponds to a fixed limit to the
number of samples taken.) the success probability indicates the number
of runs in which all 6 modes were found. When only $30\%$ of data
is available, the NERA$\ell$SFA without interpolation cannot find
all two significant modes within 200 iterations.}
\end{table}

\addtocounter{figure}{+2}

We tried another example of signal when $N=100$. NERA$\ell$SFA without
interpolation techniques fails to recover the signal with high probability
if more than $45\%$ data are unavailable. In contrast, with the help
of interpolation technique, the NERA$\ell$SFA can always recover
the signal with only $25\%$ available data.

Experiments also show that for NERA$\ell$SFA with interpolation technique,
the total number of available data, instead of the percentage of available
data determines the success probability. On the contrary, The success
of NERA$\ell$SFA without interpolation is determined by the percentage.

\subsection{Experiments to Recover Noisy Signals}

To recover a signal from very noisy data is a challenging problem.
The following tests are done for $S(t)=\sum_{i=1}^{B}c_{i}\phi_{\omega_{i}}+\nu$,
$B=6$, $\epsilon=0.02$, $N=2^{17}$, $p=0.6$, and different standard
deviation $\sigma$ for noise. The amplitude of noise is still multiplied
by a factor of $1/\sqrt{N}$. As Table \ref{tab:noise} shows, NERA$\ell$SFA
excels at extracting information from noisy data even in the case
of small signal to noise ratio.

\begin{table}
\begin{center}\begin{tabular}{|c|c|c|c|c|c|}
\hline 
$\sigma$&
 SNR&
 Time of NERA$\ell$SFA &
 Time of NERA$\ell$SFA &
Relative Error &
 Success \tabularnewline
&
 (dB)&
 (+sampling)&
 ( w/o sampling)&
 ($\%$)&
 probability \tabularnewline
\hline
0 &
 -&
 0.48&
 0.21&
 0.02&
 100\%\tabularnewline
\hline
0.5 &
 27.60&
 0.56&
 0.22&
 2.00&
 100\%\tabularnewline
\hline
1.0&
 15.56&
 0.87&
 0.32&
 4.50&
 90\%\tabularnewline
\hline
1.5&
 8.53&
 3.94&
 1.59&
 5.83&
 80\%\tabularnewline
\hline
2.0&
 3.52&
 4.78&
 1.86&
 7.67&
 50\%\tabularnewline
\hline
2.5&
 -0.35&
 7.96&
 2.14&
 8.50&
 30\% \tabularnewline
\hline
\end{tabular}\end{center}

\caption{\label{tab:noise}Experiments with fixed $B=6$, $N=2^{17}$, $p=0.6$,
and varying noise levels. For each noise level, 10 different runs
were carried out; the average result is shown. In each run, the number
of iterations is limited to 200; (this also corresponds to a fixed
limit to the number of samples taken.) the success probability indicates
the number of runs in which all 6 modes were found. The average relative
error is the error of reconstructed signal with respect to the original
signal.}
\end{table}

\addtocounter{table}{+1}

\section{Conclusion}

We provide a sublinear sampling algorithm that recovers, with high
probability, a $B$-term Fourier representation for an unevenly spaced
signal. It is faster than any existed methods for processing sparse
signals of large size. Moreover, it recovers the signal in the situation
of large percentage of missing data or small signal to noise ratio.

\section{Acknowledgments}

For many helpful suggestions and discussions, I would thank my adviser
Ingrid Daubechies. In addition, I thank Weinan E, Anna Gilbert, Martin
Strauss for their suggestions.

\section*{Appendix}

\subsection*{How to interpolate the two dimensional data to get values for missing
points}

In one dimension, values of missing points can be interpolated by
its few nearest left and right available neighbors. The idea can be
extended to higher dimensional cases with more techniques.

For instance, in two dimensions, we first find four nearest available
neighbors of a missing point in each quadrant. Suppose a missing point
is $(x,y)$, its four neighbors are $(x_{1},y_{1})$, $(x_{2},y_{2})$,
$(x_{3},y_{3})$, $(x_{4},y_{4})$. The weights of neighbors can be
derived by solving the following linear system of equations.

\begin{equation}
\left(\begin{array}{cccc}
x_{1} & x_{2} & x_{3} & x_{4}\\
y_{1} & y_{2} & y_{3} & y_{4}\\
x_{1}y_{1} & x_{2}y_{2} & x_{3}y_{3} & x_{4}y_{4}\\
1 & 1 & 1 & 1\end{array}\right)\left(\begin{array}{c}
w_{1}\\
w_{2}\\
w_{3}\\
w_{4}\end{array}\right)=\left(\begin{array}{c}
x\\
y\\
xy\\
1\end{array}\right)\label{array1}\end{equation}

However, the matrix in (\ref{array1}) could be singular. In this
case we choose the three nearest neighbors in different quadrants
and use the following equations:

\begin{equation}
\left(\begin{array}{ccc}
x_{1} & x_{2} & x_{3}\\
y_{1} & y_{2} & y_{3}\\
1 & 1 & 1\end{array}\right)\left(\begin{array}{c}
w_{1}\\
w_{2}\\
w_{3}\end{array}\right)=\left(\begin{array}{c}
x\\
y\\
1\end{array}\right)\label{array2}\end{equation}
 \label{array2}

The time to locate those nearest neighbors and compute corresponding
weights is considered a part of precomputation and excluded from total
running time.

Note that we can use geometrical arguments to simplify the pre-computation
of the weights. One easily sees that the system of equations (\ref{array1})
is translation invariant: the two linear system of equations

\[
\left(\begin{array}{cccc}
x_{1}+l & x_{2}+l & x_{3}+l & x_{4}+l\\
y_{1}+p & y_{2}+p & y_{3}+p & y_{4}+p\\
(x_{1}+l)(y_{1}+p) & (x_{2}+l)(y_{2}+p) & (x_{3}+l)(y_{3}+p) & (x_{4}+l)(y_{4}+p)\\
1 & 1 & 1 & 1\end{array}\right)\left(\begin{array}{c}
w_{1}\\
w_{2}\\
w_{3}\\
w_{4}\end{array}\right)=\left(\begin{array}{c}
l\\
p\\
lp\\
1\end{array}\right)\]
 and \[
\left(\begin{array}{cccc}
x_{1} & x_{2} & x_{3} & x_{4}\\
y_{1} & y_{2} & y_{3} & y_{4}\\
x_{1}y_{1} & x_{2}y_{2} & x_{3}y_{3} & x_{4}y_{4}\\
1 & 1 & 1 & 1\end{array}\right)\left(\begin{array}{c}
w_{1}\\
w_{2}\\
w_{3}\\
w_{4}\end{array}\right)=\left(\begin{array}{c}
0\\
0\\
0\\
1\end{array}\right)\]
 have the same solutions for any $l$ and $p$. That means the location
of the missing points does not influence the weights. Only the geometrical
shape and relative distance of the available neighbors of a missing
point matters.

Thus, we compute weights for the geometrical shapes of available neighboring
points which occur most often. As we go through every missing point,
we check if the shape of its neighboring available points matches
those popular ones; if it does, we can directly get the weights without
computation. This saves a huge amount of work, especially when $p$
is large.

\begin{figure}
\includegraphics[%
  width=7cm]{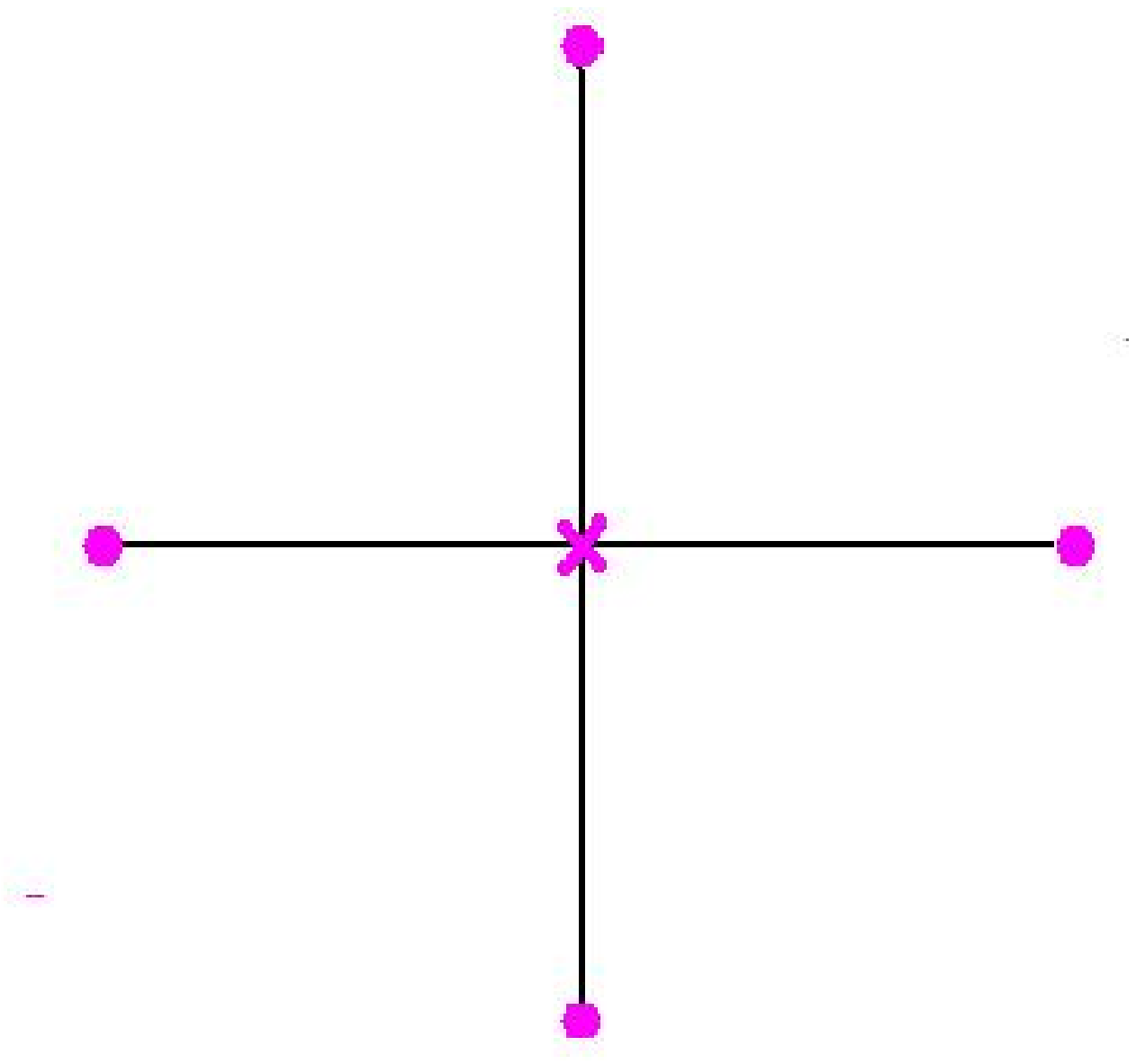}\includegraphics[%
  width=7cm]{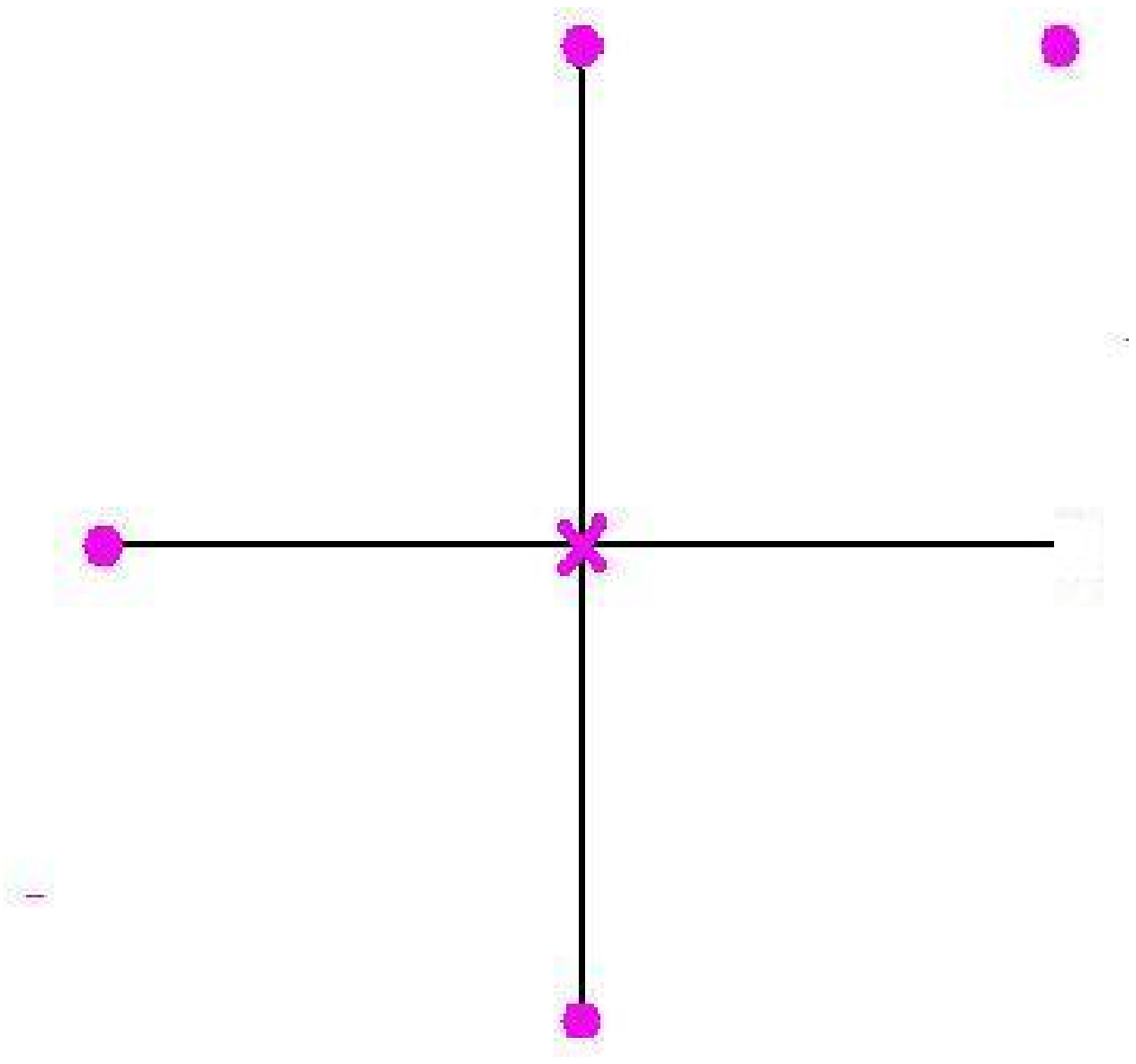}

\caption{\label{fig:shape}Some geometrical shapes of available neighboring
points that occur most often. A missing point (denoted by a small
cross) is at the center of the cross. Available points are denoted
by dots. Left: the four available neighbors are located in the shape
of cross. The distances of each neighbor to the missing point are
equal. Right: almost the same as configuration in the left side, except
one point moved off to the diagonal. }
\end{figure}

For example, if the four neighboring points are located in the shape
of a cross with the missing point as their center, as the left side
of Figure \ref{fig:shape} shows, then all of the weights are equal
to one quarter. This situation happens with probability $p^{4}$,
which is almost $2/3$ when $p=0.9$. Another often occurring case
typically has one of the four neighbors of the previous configuration
moved off to the diagonal (see the right side of Figure \ref{fig:shape}),
which happens with probability $4p^{4}(1-p)(2-p)$, i.e. about $28\%$
when $p=0.9$. In this case, the two neighbors on the same line as
the mirroring points have a weight 0.5 respectively; the other two
points have weight zero. Table \ref{tab:shape} shows the probabilities
of these two situations as $p$ varies.

\begin{table}
\begin{center}\begin{tabular}{|c|c|c|c|}
\hline 
$p$&
 $p^{4}$&
 $4p^{4}(1-p)(2-p)$&
 sum:$p^{4}+4p^{4}(1-p)(2-p)$\tabularnewline
\hline
1&
 $100\%$&
 0&
 $100\%$\tabularnewline
\hline
0.9&
 $65\%$&
 $29\%$&
 $94\%$\tabularnewline
\hline
0.8&
 $41\%$&
 $39\%$&
 $80\%$\tabularnewline
\hline
0.7&
 $24\%$&
 $37\%$&
 $61\%$\tabularnewline
\hline
0.6&
 $13\%$&
 $29\%$&
 $42\%$\tabularnewline
\hline
0.5&
 $6\%$&
 $19\%$&
 $25\%$ \tabularnewline
\hline
\end{tabular}\end{center}

\caption{\label{tab:shape}Two possibilities corresponding to the geometrical
shapes in Figure 9. The parameter $p$ is the percentage of available
data. The left side of Figure 9 happens with probability $p^{4}$;
the right side appears with probability $4p^{4}(1-p)(2-p)$.}
\end{table}

\end{document}